\documentclass[12pt]{amsart}
\usepackage{amsmath,
amsfonts,
amsthm,
amssymb,
mathrsfs}

\newcommand{\dom}{\operatorname{dom}}

\newcommand{\comment}[1]{}

\newcommand{\Pbb}{\mathbb{P}}
\newcommand{\rest}{\upharpoonright}
\newcommand{\Ht}{\operatorname{ht}}

\newcommand{\Seq}[1]{\langle #1 \rangle}

\newcommand{\ZFC}{\mathrm{ZFC}}
\newcommand{\GCH}{\mathrm{GCH}}
\newcommand{\CH}{\mathrm{CH}}

\theoremstyle{plain}
\newtheorem{thm}{Theorem}[section]
\newtheorem{lem}[thm]{Lemma}
\newtheorem{prop}[thm]{Proposition}

\newtheorem{claim}[thm]{Claim}

\theoremstyle{definition}
\newtheorem{defn}[thm]{Definition}

\begin{document}

\author[H. Lamei Ramandi]{Hossein Lamei Ramandi}

\address{}

\title[Rigidity of Souslin trees and generic
branches ]{On the rigidity of Souslin trees and
their generic branches}

\subjclass[]{}

\email{{\tt hlamaira@exchange.wwu.de}}
 \begin{abstract}
We show it is consistent that there is a Souslin
tree $S$ such that after forcing with $S$, $S$ is Kurepa and
for all clubs $C \subset \omega_1$, $S\rest C$ is
rigid. This answers the questions in \cite{Fuchs}.
Moreover, we show it is consistent with
$\diamondsuit$ that for every Souslin tree $T$
there is a dense
$X \subseteq T$ which does not contain a copy of
$T$. This is related to a question due to
Baumgartner in \cite{baumgartner_Other}.
\end{abstract}
 
 \maketitle
\section{Introduction}

Recall that an $\omega_1$-tree is said to be
Souslin if it has no uncountable chain or
antichain.
In \cite{Fuchs} and \cite{FuchsHamkins}, Fuchs and
Hamkins considered various notions of rigidity of
Souslin trees and  
studied the following question: How many generic
branches can Souslin trees introduce, when they satisfy certain
rigidity requirements?
In \cite{Fuchs}, Fuchs asks a few questions which motivate the following theorem.
\begin{thm}\label{main}
It is consistent with $\GCH$ that there is a
Souslin tree $S$ such that
$\Vdash_S ``S$ is Kurepa and $S \rest C$ is rigid
for every club $C \subset \omega_1$".

\end{thm} 
\noindent
Theorem \ref{main} answers all questions in
\cite{Fuchs}.
We refer the reader to \cite{Fuchs} and
\cite{FuchsHamkins} for motivation and history.

In \cite{baumgartner_Other}, Baumgartner proves
that under $\diamondsuit^+$ there is a
lexicographically ordered Souslin tree which is
minimal as a tree and as an uncountable linear
order.
At the end of his construction he asks the
following question: Does there exist a minimal
Aronszajn line if $\diamondsuit$ holds?
This question is not settled here but motivates
the following proposition.
\begin{prop} \label{prop}
It is consistent with $\diamondsuit$ 
that if $S$ is a Souslin tree then there is a
dense $X \subset S$
which does not contain a copy of $S$.
\end{prop}    
\noindent
Proposition \ref{prop} shows it is impossible to
follow the same strategy as Baumgartner's in
\cite{baumgartner_Other},
in order to show $\diamondsuit$ implies that there
is a minimal Aronszajn line.
More precisely, it is impossible to find a
lexicographically ordered Souslin tree which is
minimal
as a tree and as an uncountable linear order. 

This paper is organized as follows. In the next section we prove Proposition 
\ref{prop}.
In the third section we introduce a Souslin tree which makes itself a Kurepa tree.
This tree is used in the last section, where we prove Theorem \ref{main}.

Let's fix some definitions, notations and conventions.
Assume $T,S$ are trees and $f : T \longrightarrow
S$ is injective.
Then $f$ is said to be an \emph{embedding} when  $t
<_T s \Longleftrightarrow f(t) <_S f(s)$.
$T$ is called an \emph{$\omega_1$-tree} if its levels are
countable and $\Ht(T)= \omega_1$.
$T$ is said to be pruned if for all $t \in T$ and $\alpha \in \omega_1 \setminus \Ht(t)$ there is $s \geq t$ such that $\Ht(s) = \alpha$.
If $t \in T$ and $\alpha \leq \Ht(t)$, $t \rest \alpha$ refers to the 
$\leq_T$ predecessor of $t$ in level $\alpha$.
$C \subset T$ is called a \emph{chain} if it consists of pairwise comparable 
elements.
A chain $b \subset T$ is called a \emph{branch} if it intersects all levels of $T$.
An $\omega_1$-tree $U$ is called \emph{minimal} if
for every uncountable $X \subset U$,
$U$ embeds into $X$. 
If $T$ is a tree and $\alpha$ is an ordinal, $T(\alpha)= \{ t \in T: \Ht(t) = \alpha
\}$ and
$T(< \alpha) = \{ t \in T: \Ht(t) < \alpha \}$.
If $A$ is a set of ordinals, $T \rest A = \{t \in T:
\Ht(t) \in A \}$.
If $t \in T$ and $U \subset T$ then $U_t = \{ u \in U: t \leq_T u \}$.
Assume $Q$ is a poset and $\theta$ is a regular cardinal. 
We say $M \prec H_\theta$ is suitable for $Q$ if $Q$ and the power set of
the transitive closure of $Q$ are in $M$.
 
\section{Minimality of Souslin trees and $\diamondsuit$}

This section is devoted to the proof of Proposition \ref{prop}.
We will use the following terminology and notation in this section.
By $N$ we mean the set of all countable infinite successor ordinals, and
$\Pbb$ refers to the countable support iteration 
$\Seq{P_i, \dot{Q}_j : i \leq \omega_2, j <
\omega_2}$, where $Q_j=2^{<\omega_1}$
for each $j \in \omega_2$.

\begin{lem} \label{nonchain}
Assume $U= (\omega_1, <)$ is a Souslin tree, $p
\in \Pbb$,
$\dot{X}$ is the canonical ${\Pbb}_1$-name for the
generic subset of $\omega_1$,
$p \Vdash ``\dot{f}$ is an embedding from $U$ to
$\dot{X}$" and for every $t \in U$ define 
$\varphi(p,t) = \{ s \in U : \exists \bar{p} \leq
p \textrm{ } \bar{p} \Vdash \dot{f}(t) = s \}.$
Then there is an $\alpha \in \omega_1$ such that
for all $t \in U \setminus U(< \alpha) $,
$\varphi(p,t)$ is not a chain.
\end{lem}
\begin{proof}
Let $Y_p = \{ y \in U : \varphi(p,y) \textrm{ is a chain} \}.$
$Y_p$ is downward closed and if it is countable we are done.
Fix $p \in \Pbb$ and assume for a contradiction
that $Y_p$ is uncountable.
Let $A_p = \{ t \in U : p \Vdash t \in \dot{X}
\textrm{ or } p\Vdash t \notin \dot{X} \}.$
$A_p$ is countable. Fix $\alpha > \sup \{ \Ht(a):
a \in A_p \}$ and
$y \in Y_p \setminus U(\leq \alpha)$.
Since $U$ is an Aronszajn tree and $\varphi(p,y) $
is a chain,
we can choose $\beta \in \omega_1 \setminus \sup
\{ \Ht(s): s \in \varphi(p,y) \}$.
For all $s \in \varphi(p,y)$, $\alpha < \Ht(s) <
\beta$
since $\emptyset \Vdash \Ht(y) \leq \Ht(\dot{f}(y))$.
Then we can extend $p$ to $q$ such that $q \Vdash
\dot{X} \cap (U(\leq \beta) \setminus U(< \alpha))
= \emptyset$, which
contradicts $p \Vdash \dot{f}(y) \in
\varphi(p,y)$.
\end{proof}

\begin{lem}\label{Silver1}
Assume $U \in \textsc{V}$ is a pruned Souslin tree and
$G \subset \Pbb$ is $\textsc{V}$-generic. Then in $\textsc{V}{[G]}$, there
is
a dense $X \subset U$ which does not have a copy
of $U$.
\end{lem}
\begin{proof}
Let $\dot{X}$ be as in Lemma \ref{nonchain}.
Since $U$ is pruned, $1_\Pbb \Vdash \dot{X} \subset U$ is
dense.
We will show $1_\Pbb \Vdash \dot{X}$ has no copy
of $U$.
Assume for a contradiction that $p \Vdash_\Pbb
\dot{f}$ is an embedding from $U$ to $\dot{X}$.
Fix a regular cardinal $\theta$ and a countable $M
\prec H_\theta$ which contains $U, p, \dot{f},
2^\Pbb$.
Also let $\Seq{D_n : n \in \omega}$ be an
enumeration of all dense open subsets of $\Pbb$ in $M$, $\delta
= M \cap \omega_1$ and $t \in U(\delta)$.
For each $\sigma \in 2^{< \omega}$, find 
$p_\sigma \in D_{|\sigma|} \cap M$, $s_\sigma$ and
$t_{|\sigma|}< t$,
such that:
\begin{enumerate}
\item if $\sigma \sqsubset \tau$ then $p_\tau \leq
p_\sigma$ and $s_\sigma \leq s_\tau$,
\item if $\sigma \perp \tau$ then $s_\sigma \perp
s_\tau$,
\item $p_\sigma \Vdash \dot{f}(t_{|\sigma|}) =
s_\sigma$.
\end{enumerate}

In order to see how these sequences are
constructed, let $t_0 < t$ be arbitrary and
$p_{\emptyset}, s_\emptyset$ be
such that $p_\emptyset \Vdash`` \dot{f}(t_0)=
s_\emptyset"$ and $p_\emptyset \in D_0 \cap M$.
Assuming these sequences are given for all $\sigma
\in 2^n$, use Lemma \ref{nonchain} to find $
t_{n+1} < t$
such that $\varphi(p_\sigma, t_{n+1})$ is not a
chain, for all $\sigma \in 2^n$.
Let $s_{\sigma \frown 0}, s_{\sigma \frown 1}$ be
in $\varphi(p_\sigma , t_{n+1}) \cap M$ such that
$s_{\sigma \frown 0} \perp s_{\sigma \frown 1}$.
Now find $p_{\sigma \frown 0} , p_{\sigma \frown
1}$ in $M \cap D_{n+1}$
which are extensions of $p_\sigma$ such that 
$p_{\sigma \frown i} \Vdash ``\dot{f}(t_{n + 1}) =
s_{\sigma \frown i}", \textrm{ for } i=0,1.$

For each $r \in 2^\omega$, let $p_r$ be a lower
bound for $\{ p_\sigma : \sigma \sqsubset r \}$
and let $b_r \subset U \cap M$ be a downward
closed chain such that
$p_r \Vdash \dot{f}[\{ s \in U : s < t \}] \subset
b_r$.
Note that $b_r$ intersects all the levels of $U$
below $\delta$.
It is obvious that $p_r$ is an $(M,\Pbb)$-generic
condition below $p$.
Moreover, if $r,r'$ are two distinct real numbers then $b_r \neq b_{r'}$.
Let $r \in
2^\omega$ such
that $U$ has no element on top of $b_r$.
Then $p_r$ forces that $\dot{f}(t)$ is not
defined, which is a contradiction.
\end{proof}

Now we are ready for the proof of Proposition
\ref{prop}.
Let $\textsc{V}$ be a model of $\ZFC + \GCH$ and
$G \subset \Pbb$ be $\textsc{V}$-generic.
Since $\Pbb$ is a countable support iteration of 
$\sigma$-closed posets of size $\aleph_1$, 
it preserves all cardinals.
The same argument as in Theorem 8.3 in
\cite{set_theory:Kunen}
shows that $\diamondsuit$ holds in
$\textsc{V}[G]$.

Let $U$ be a Souslin tree in $\textsc{V}[G]$. 
For some $\alpha \in \omega_2$, $U \in
\textsc{V}[G \cap P_\alpha]$ since $|U| =
\aleph_1$.
Let  $\dot{R}$ be the canonical $P_\alpha$-name 
such that $\Pbb = P_\alpha * \dot{R}$.
Then $1_{P_\alpha} \Vdash \dot{R}$ is isomorphic
to $\Pbb$.
By Lemma \ref{Silver1}, there is a dense $X
\subset U$ in $\textsc{V}$[G]
which has no copy of $U$, as desired.

\section{A Souslin tree with many generic branches
}
\begin{defn}
The poset $Q$ is the set of all $p=(T^p, \Pi_p)$
such that:
\begin{enumerate}
\item $\Delta_p \in \omega_1$ and $T^p = (\Delta_p, \leq_p)$
is a countable binary tree of height $\alpha_p$ such that
for all $t \in T^p$ and for all $\beta \in
\alpha_p \setminus \Ht_{T^p}(t)$ there is $s \in
T^p(\beta)$
with $t<_{T^p} s$.
\item $\Pi_p = \Seq{\pi_\xi^p: \xi \in D_p}$ where
$D_p \subset \omega_2$ is countable and
for each $\xi \in D_p$ there are $x,y$ of the same
height in $T^p$ such that $\pi_\xi^p : (T^p)_x
\longrightarrow (T^p)_y$
is a tree isomorphism.
\end{enumerate}
We let $q \leq p $ if $T^q$ end-extends $T^p$,
$D_p \subset D_q$ and
for all $\xi \in D_p$, $\pi_\xi^q \rest T^p = \pi_\xi^p $.

\end{defn}

\begin{lem}
$Q$ is $\sigma$-closed. Moreover if $\CH$ holds,
$Q$ has the $\aleph_2$-cc.
\end{lem}
\begin{proof}
The first part of the lemma is obvious.
Assume $A \in Q^{\aleph_2}$.
By thinning $A$ out, we can assume that for all $p,q$ in $A$, $T^p =T^q$,
$\{ D_p : p \in A \}$ is a $\Delta$-system with root $R$ and
$|\{ \Seq{\pi_\xi^p : \xi \in R}: p \in A \}|=1$. 
Now all $p,q$ in $A$ are compatible.
\end{proof}

\begin{lem}\label{Souslin}
If $T = \bigcup_{p \in G}T^p$ for a generic $G \subset Q$, then $T$ is Souslin.
\end{lem}
\begin{proof}
Obviously $T$ is an $\omega_1$-tree.
Let $\tau$ be a $Q$-name and $p \Vdash_Q ``\tau
\subset T$ is a maximal antichain\rq\rq{}.
We show $p \Vdash \tau$ is countable.
Let $M \prec H_\theta$ be countable, $\theta$
regular and $2^Q, \tau$ be in $M$.
Let $\Seq{p_n = (T_n , \Pi_n) :n \in \omega}$, be
a descending $(M,Q)$-generic sequence with $p_0 = p$.
Let $\pi_\xi^{p_n}= \pi_\xi^n$, $\delta = M \cap \omega_1$,
and $R = \bigcup_{n \in \omega} T_n$. 
So $\Ht(R) = \delta$ and 
$M \cap \omega_2 = \bigcup_{n \in \omega}D_{p_n}$.
Let $\mathcal{F}$ be the
set of all finite compositions of functions of the form
$\bigcup_{n \in \omega} \pi_\xi^n$ with $\xi \in M \cap \omega_2 $.
Let $\Seq{f_n : n \in \omega}$ be an enumeration of $\mathcal{F}$
with infinite repetition and $A= \{ t \in R: \exists n \in \omega$
$(p_n \Vdash t \in \tau) \}$.
Observe that for all $t \in R$ there is $a \in A$ such that $a,t$ are comparable.

Let $\Seq{\alpha_m: m \in \omega}$ be an
increasing cofinal sequence in $\delta$.
For each $t \in R$ we build an increasing sequence $\bar{t} = \Seq{t_m : m \in \omega }$ as follows.
Let $t_0 =t$.
Assume $t_m$ is given. 
If $R_{t_m} \cap \dom(f_m) = \emptyset$, choose $t_{m+1} > t_m$
with $\Ht(t_{m +1}) > \alpha_m$.
If $R_{t_m} \cap \dom(f_m) \neq  \emptyset$, let $s \in \dom(f_m) \cap R_{t_m}$.
Let $a \in A$ such that $a, f_m(s)$ are comparable.
Let $x = \max \{f_m(s),a \}$ and $t_{m+1}> f^{-1}_m(x)$ with
$\Ht(t_{m+1})> \alpha_m$.
Let $b_t$ be the downward closure of $\bar{t}$.

Let $B = \{ f_n[b_t] : t \in R \textrm{ and } n \in \omega \}$.
Let $q$ be the lower bound for $\Seq{p_n: n \in \omega}$ described as follows.
$T^q = R \cup T^q(\delta)$ and for each cofinal 
branch $c \subset R$ there is a unique $y \in T^q(\delta)$ above $c$ if and only if $ c \in B $.
For each $\xi \in M \cap \omega_2$, 
let $\pi_\xi^q \rest R = \bigcup_{n \in \omega}
\pi_\xi^n$.
Note that this determines $\pi_\xi^q$ on
$T^q(\delta)$ as well and
$\pi_\xi^q (y)$ is defined for all $y \in
T^q(\delta)$.

The condition $q$ forces that for each $y \in
T(\delta) = T^q(\delta)$
there is $a \in A$
with $a < y$.
In other words $q$ forces that $\tau = A$.
Since $p$ was arbitrary, 
$1_Q$ forces that every maximal antichain has to
be countable.
\end{proof}

From now on $T$ is the same tree as in Lemma
\ref{Souslin}.
For each $\xi \in \omega_2$ let 
$\pi_\xi = \bigcup_{p \in G} \pi^p_\xi$, where $G
\subset Q$ is generic.
Observe that if $x \in \dom(\pi_\xi) \cap \dom(\pi_\eta)$
and $\xi \neq \eta$ are ordinals
then there is $\alpha > \Ht(x)$ such that 
for all $y \in T(\alpha) \cap T_x$, $\pi_\xi (y)
\neq \pi_\eta(y)$.
So forcing with $T$ makes $T$ Kurepa.

\section{Highly rigid dense subsets of $T$}
In this section we show the tree $T$, in the
forcing extensions by $
P = (2^{<\omega_1}, \supset)$, has dense subsets which are witnesses for
Theorem
\ref{main}.

\begin{lem}\label{unused}
Let $U=(\omega_1, <)$ be a pruned
Souslin tree and
$S \subset \omega_1$ be generic for $P$.
Then in $\textsc{V}[S]$ the following hold.
\begin{enumerate}
\item $S$ is a Souslin tree when it is considered with the inherited order from 
$U$.
\item $S \subset U$ is dense.
\item For all clubs $C \subset \omega_1$, $S \rest C$ is
rigid.
\end{enumerate}

\end{lem}
\begin{proof}
In order to see that $S$ is Souslin, note that $\sigma$-closed posets do not 
add uncountable antichains to Souslin trees. Moreover by standard 
density arguments $S \subset U$ is dense.

Assume for a contradiction $p \Vdash_P ``\dot{f}:
\dot{S} \rest \dot{C}
\longrightarrow \dot{S} \rest \dot{C}$ is a
nontrivial tree embedding.\rq\rq{}
Let $\Seq{M_\xi : \xi \in \omega +1}$ be a
continuous $\in$-chain of
countable elementary submodels of $H_\theta$
where $\theta$ is regular and $p, \dot{f}, 2^U $ are in $ M_0$.
For each $\xi \leq \omega$, let $\delta_\xi =
M_\xi \cap \omega_1$ and
$t \in U({\delta_\omega})$.
Let $t_n = t \rest \delta_n$.
For each $\sigma \in 2^{< \omega}$ we find
$q_\sigma \in M_{|\sigma|+1} \cap P$, $s_\sigma$ such
that:
\begin{enumerate}
\item $q_0 \leq p$, and if $\sigma \subset \tau$
then $q_\tau \leq q_\sigma$,
\item $q_\sigma$ is $(M_{|\sigma|}, P)$-generic
and
$q_\sigma \subset M_{|\sigma|}$,
\item $q_\sigma$ forces that
$\dot{f}(t_{|\sigma|-1}) = s_\sigma$,
\item if $\sigma \perp \tau$ then $s_\sigma \perp
s_\tau$,
\item if $\sigma \subset \tau$ then $q_\tau$
forces that
$t_{|\sigma|} \in \dot{S} \rest \dot{C}$.
\end{enumerate}
Assuming $q_\sigma$ and $s_\sigma$ are given for
all $\sigma \in 2^{n}$,
we find $q_{\sigma \frown 0}, q_{\sigma \frown 1},s_{\sigma \frown 0},$
and $
s_{\sigma \frown 1}$.
Let $\bar{q}_\sigma = q_\sigma \cup \{( t_n,1)
\}$.
Obviously, $\bar{q}_\sigma \Vdash t_n \in \dot{S}
\rest \dot{C}$ and for all $\sigma \in 2^n$,
$\{ s \in U : \exists r \leq \bar{q}_\sigma$ $r
\Vdash \dot{f}(t_n)=s\}$ is uncountable.
In $M_{n+1}$, find $r_0,r_1$   
below $\bar{q}_\sigma$ and $s_{\sigma \frown 0},
s_{\sigma \frown 1}$
such that $s_{\sigma \frown 0} \perp s_{\sigma \frown 1}$
and $r_i \Vdash ``\dot{f}(t_n) = s_{\sigma \frown i}$.\rq\rq{}
Let $q_{\sigma \frown i} < r_i$ be $(M_{n+1}, P)$-generic
with $q_{\sigma \frown i} \subset M_{n +1}$, and $q_{\sigma \frown i} \in M_{n +2}$.

Let $r \in 2^\omega$ such that $\{ s_\sigma : \sigma \subset r \}$
does not have an upper bound in $U$.
Let $p_r$ be a lower bound for $\{p_\sigma : \sigma \subset r \}$.
Then $p_r$ forces that $\dot{f}(t)$ is not defined which is a contradiction.
\end{proof}

\begin{lem}\label{min_extension}
Suppose $M$ is suitable for $Q$ and
$\delta = M \cap \omega_1$.
Let $\Seq{q_n : n \in \omega}$ be a decreasing $(M,Q)$-generic sequence.
Define a condition $q \in Q$ by setting 
$T^q = \bigcup_{n \in \omega} T^{q_n}$, 
$D_q = \bigcup_{n \in \omega}D_{q_n}$
and for each $\xi \in D_q$
let $\pi_\xi^q = \bigcup_{n \in \omega}\pi_\xi^{q_n}$.
Also let $\Pi_q = \Seq{\pi_\xi^q: \xi \in D_q}$.
Let $\mathcal{F}$ be the set of all finite compositions of functions of the form 
$\pi^q_\xi$ with $\xi \in D_q$. 
Assume $m \in \omega$ and $\Seq{b_i : i \in m}$ are branches through $T^q$.
Then there is an extension $q\rq{} \leq q$ such that 
$\alpha_{q\rq{}} \geq \delta +1$ and for all branches $c \subset T^q$, 
$c$ has an upper bound  iff
for some $f \in \mathcal{F}\textrm{ and } i \in m,$ $ f(b_i)$ is cofinal in  $c.$
\end{lem}
\begin{proof}
Note that $D_q= M \cap \omega_2$ and $\alpha_q = \delta$.
Let $T^{q\rq{}} \rest \delta = T^q$. 
Let $B= \{ f(b_i) : i \in m \textrm{ and } f \in \mathcal{F} \}$.
Obviously $B$ is countable and we can fix an enumeration of $B$ with $n \in \omega$.
Let $T^{q\rq{}}(\delta +1) = [ \delta, \delta + \omega)$
and put $\delta + n$ on top of the $n$\rq{}th element in $B$.
It is obvious how we should extend $\Pi_q$ to $\Pi_{q\rq{}}$
with $D_q = D_{q\rq{}}$.
\end{proof}

\begin{lem}\label{nonchain2}
Let $G \subset Q$ be $\textsc{V}$-generic, $p \in P$ and
$\dot{S}$ be the canonical $P$-name for the generic subset 
of $\omega_1$.
Let $\dot{f}, \dot{C}$ be $P * T$-names in $\textsc{V}[G]$ and $t,x,y$
be pairwise incompatible in $T$.
Suppose $(p,t)$ forces $\dot{f}$ is an embedding from 
$\dot{S}_x \rest \dot{C}$ to $\dot{S}_y \rest \dot{C}$.
For every $u \in T_x$ define 
$\psi(p,t,u) = \{s \in T : \exists t\rq{} > t \textrm{ }
\exists \bar{p} \leq p \textrm{ }
(\bar{p}, t\rq{}) \Vdash [u \in \dot{S}_x \rest \dot{C} \wedge
\dot{f}(u) = s] \}.$
Then for any $u \in T_x$ there is $u\rq{} > u$ such that 
$\psi(p,t, u\rq{})$ is not a chain.
\end{lem}
\begin{proof}
Fix $p,t,u$ as above and assume for a contradiction that for all $u\rq{} > u$ in $T$,
$\psi(p,t,u\rq{})$ is a chain. 
Since $T$ is ccc, without loss of generality we can assume that 
for all $q \in P$ and $ \alpha \in \omega_1$, there is $\bar{q} \leq q$
such that $(\bar{q},1_T)$ decides the statement $\alpha \in \dot{C}$.
For each $q \in P,r \in T, v \in T$ let
$\alpha_{q,r,v} = \sup \{ \Ht_T(s) : s \in \psi(q,r, v) \}$.
Note that if $\bar{q} \leq q$ and $\bar{r} \geq r$ then $\psi(\bar{q},\bar{r},v) \subseteq \psi(q,r,v)$
and $\alpha_{\bar{q},\bar{r},v} \leq \alpha_{q,r,v}$.

Let $M_0, M_1$ be countable elementary submodels of 
$H_\theta$, $\theta$ be a regular cardinal and 
$\{p,t,u,x,y, \dot{f}, \dot{C}\} \in M_0 \in M_1$.
Suppose $\Seq{p_n: n \in \omega}$ is an $(M_0,P)$-generic sequence 
which is in $M_1$ and $p_0 \leq p$.
Let $p\rq{} = \bigcup_{n \in \omega} p_n$ and 
$\delta_i = M_i \cap \omega_1$, for $i \in 2$.
Note that $p\rq{} \Vdash \delta_0 \in \dot{C}$.

Let $\bar{p} < p\rq{}$ such that:
\begin{enumerate}
\item $\bar{p} \Vdash \forall v \in T_x \cap (M_1 \setminus M_0) \textrm{ } [v \in \dot{S}]$
\item $\bar{p} \Vdash \forall v \in T_y \cap (M_1 \setminus M_0) \textrm{ }
[v \notin \dot{S}]$.
\end{enumerate}

Let $u_0 > u$ be in $T(\delta_0)$.
Since $\bar{p}$ is $(M_0,Q)$-generic, it forces that $\delta_0 \in \dot{C}
\wedge u_0 \in \dot{S} \wedge \Ht_{\dot{S}}(u_0) = \delta_0$.
In particular, by elementarity of $M_0$ and basic facts on ordinal arithmetic,
 $\bar{p} \Vdash u_0 \in \dot{S}_x \rest \dot{C}$.

Suppose $q< \bar{p}, r > t$ such that $(q,r)$ decides $\dot{f}(u_0)$.
Then the condition $(q,r)$ forces that $\Ht(\dot{f}(u_0)) \geq \delta_1$.
So,
$\delta_1 \leq \alpha_{\bar{p},t,u_0 } \leq \alpha_{p\rq{}, t, u_0} \in M_1.$
But this is a contradiction.
\end{proof}

In the next lemma we use the following standard fact: If $U$ is a Souslin tree and $X \subset U$ is uncountable and downward closed, then there is $x \in U$
such that $U_x \subset X$. In order to see this 
assume for all $v \in U,$ $U_v$ is not contained in $X$.
Let $A$ be the set of all minimal $a$ outside of $X$. 
Observe that $A$ is an uncountable antichain, contradicting the fact that 
$U$ was Souslin.
Lemma \ref{last} finishes the proof of Theorem
\ref{main}.

\begin{lem} \label{last}
Assume $G*S*b$ is $\textsc{V}$-generic for $Q*P*\dot{T}$.  
Let $x, y$ be incomparable in $T$.
Then in $\textsc{V}[G*S*b]$ for all clubs $C \subset
\omega_1$, $S_x \rest C$ does not embed into $S_y \rest C$.
\end{lem}
\begin{proof}
Assume for a contradiction that $(q_0,p_0,t_0)$ is a condition in 
$Q*P *\dot{T}$ which forces 
$\dot{f}: \dot{S}_x \rest \dot{C}\longrightarrow \dot{S}_y \rest \dot{C}$
is a tree embedding and $x,y$ are incompatible in $T$.
Note that $\dot{f}(\dot{S}_x)$ is an uncountable subset of $\dot{T}_y$
and $\dot{T}$ is a Souslin tree in $\textsc{V}[G][S]$.
So the downward closure of $\dot{f}(\dot{S}_x)$ contains  $\dot{T}_z$
for some $z > y$. 
Therefore, by extending $x,y,(q_0,p_0,t_0)$ if necessary,
we can assume that $\dot{f}(\dot{S}_x)$ is dense in $\dot{S}_y$.

Again by extending $x,y,(q_0,p_0,t_0)$ we may assume  
$(q_0, p_0,t_0)$ $\Vdash$ $[x,y$ are in $\dot{S} \rest \dot{C}$ and $\dot{f}(x) = y]$.
Furthermore, by extending $t_0$ if necessary we can assume that 
$\Ht(t_0) > \Ht (y)$ and $x,y,t_0$ are pairwise incomparable.
Since $T$ is a ccc poset we can assume that for
all $\alpha \in \omega_1$, for all $u,v$ in $T$ and for all $(a,b) \in P*Q$
we have $(a,b,u) \Vdash \alpha \in \dot{C}$
$\longleftrightarrow$ $(a,b,v) \Vdash \alpha \in \dot{C}$.

Let $M$ be a countable elementary submodel of $H_\theta$
such that $\theta$ is regular and $(q_0,p_0,t_0), \dot{f}$ are in $ M$.
Let $\Seq{q_n : n \in \omega}$ be a decreasing $(M,Q)$-generic sequence.
Define $q \in Q$ as in Lemma \ref{min_extension}.
Let $\mathcal{F}$ be the set of all finite compositions of 
functions of the form $\pi^q_\xi$ with  $\xi \in M \cap \omega_2$.
Let $\Pi_q=\Seq{\pi_\xi^q: \xi \in M \cap \omega_2}$.
Obviously, $q$ is an $(M,Q)$-generic condition.
Let $\Seq{g_n : n \in \omega}$ be an enumeration of $\mathcal{F}$ with infinite repetition.
Let $\Seq{\gamma_n : n \in \omega}$ be an increasing cofinal sequence in 
$\delta = M \cap \omega_1$ with $\gamma_0 = 0$. 
                           
We find a decreasing sequence $\Seq{p_n \in P \cap M: n \in \omega }$ and increasing sequences   
$\Seq{\delta_n \in  \delta : n \in \omega}$,
$\Seq{t_n \in T^q: n \in \omega}$, $\Seq{u_n \in T^q: n \in \omega}$
$\Seq{s_n \in T^q : n \in \omega }$  such that:
\begin{enumerate}
\item $\delta_n \geq \gamma_n$ for all $n \in \omega,$ 
\item $(q,p_n.t_n) \Vdash \min \{ \Ht_{\dot{S}}(s_n), \Ht_{\dot{S}}(u_n), 
\dom (p_n) \} \geq \delta_n$,
\item  $\Ht_{T^q}(t_n) \geq \Ht_{T^q}(s_n) +1$,
\item $(q,p_n,1_{T^q}) \Vdash \delta_n \in \dot{C}$,
\item \label{forcing_relation} $(q, p_n,t_n)$   $\Vdash$  $\dot{f}(u_n) = s_n$,
\item if $n \in \omega \setminus 1$ and 
$t_{n -1} \in \dom(g_n)$ then $g_{n} (t_n) \perp s_{n}$,
\item if $n \in \omega \setminus 1 $ and $u_{n-1} \in \dom(g_n)$ then $g_n(u_n) \perp s_n$.
\end{enumerate}

We let $u_0 =x, s_0=y, \delta_0 \in \omega_1$ such that 
$(q,p_0, t_0)$ forces that  $\min \{\Ht_{\dot{S}}(x), \Ht_{\dot{S}}(y), \alpha_{p_n} \} = \delta_0$. 
It is easy to see that this choice together with $p_0,t_0$ will satisfy the  corresponding conditions.
For given $p_n, t_n, s_n, u_n, \delta_n$ 
we introduce $p_{n+1}, t_{n+1}, s_{n+1}, u_{n + 1}, \delta_{n+1}$.

If $t_n \notin \dom(g_{n + 1}) $ let $v = s_n$.
If $t_n \in \dom(g_{n + 1}) $,  let $v\geq s_{n}$ such that $v \perp g_{n+1}(t_n)$.
Such a $v$ exists because $\Ht(t_n) > \Ht(s_n)$, $g_{n +1}$ is level preserving and the tree $T^q$ is binary.
\begin{claim}
There are $t_n\rq{} > t_n$, $p_n\rq{} < p_n ,$ $u_n\rq{} > u_n$
such that  if $u_n \in \dom(g_{n+1})$ then
$(q,p_n\rq{}, t_n\rq{})$ forces 
$[u_n\rq{} \in \dom(\dot{f})$  $\wedge$ $v < \dot{f}(u_n\rq{})$  $\wedge$
$\dot{f}(u_n\rq{}) \perp g_{n+1}(u_n\rq{})]$.
\end{claim}
\begin{proof}[Proof of Claim]
Assume $u_n \in \dom(g_{n+1})$.
Recall that $\dot{f}(\dot{S}_x)$ is forced to be dense in $\dot{S}_y$.
Let $\bar{p}_n \leq p_n , \bar{t}_n \geq t_n , a_0 > u_n, v\rq{} > v$
such that $(q, \bar{p}_n, \bar{t}_n) \Vdash \dot{f}(a_0) = v\rq{}$.
This is possible because $q$ is $(M,Q)$-generic. 
Let $a > a_0$, $t_n^0, t_n^1$ be extensions of 
$\bar{t}_n$, and $p_n^0,p_n^1$ be extensions of $\bar{p}_n$
such that $(q,p_n^i, t_n^i) \Vdash \dot{f}(a) = s_n^i$ where $i \in 2$ and $s_n^0 \perp s_n^1$.
Again, this is possible because of Lemma \ref{nonchain2} and the fact that $q$ is $(M,Q)$-generic.
Let $a\rq{} > a$ such that $\Ht(a\rq{}) > \max \{ \Ht(s_n^0), \Ht(s_n^1) \}$.
Fix $i \in 2$ such that $g_{n+1}(a\rq{}) \perp s_n^i$.
Then for all $e > a\rq{}$, $(q,p_n^i, t_n^i)$ forces that
if $e \in \dom(\dot{f})$ then $\dot{f}(e) > s_n^i$. 
Moreover it forces that $g_{n+1}(e) \perp  s_\sigma^{i}$.
Therefore, $(q,p_n^i, t_n^i)\Vdash [\forall e>a\rq{} \textrm{ } e \in \dom(\dot{f})
\longrightarrow g_{n+1}(e) \perp \dot{f}(e)]$.
Let $u_n\rq{} > a\rq{}$, $p_n\rq{} <p_n^i$ and $t_n\rq{} > t_n^i$
such that $(q,p_n\rq{},t_n\rq{}) \Vdash [u_n\rq{} \in \dom(\dot{f})]$.
Then this condition will also force  
$\dot{f}(u_n\rq{}) \perp g_{n+1}(u_n\rq{})$ and $ v < \dot{f}(u_n\rq{})$.
\end{proof}

Fix $p_n\rq, t_n\rq{}, u_n\rq{}$ as in the claim above.
By extending $p_n\rq{}$ if necessary, we can assume that $(q,p_n\rq{},1_{T^q})$
decides the $\gamma_{n+1}$\rq{}st element of $\dot{C} \setminus \delta_n$ and we let $\delta_{n+1}$ 
be this ordinal. Let $u_{n+1} > u_n\rq{}$ such that for some $p_{n+1}< p_n\rq{}$ with 
$\dom(p_{n+1}) \geq \delta_{n+1}$,  the condition $(q,p_{n+1},1_{T^q})$ forces that 
$u_{n+1} \in \dot{S} \rest \dot{C}$ and $ \Ht_{\dot{S}}(u_{n+1}) \geq \delta_{n+1}$.
Let $r > t_n\rq{}$.
By extending $(q,p_{n+1},r)$ if necessary, we can assume this condition decides $\dot{f}(u_{n+1})$.
Let $s_{n+1} \in T^q$ such that $(q,p_{n+1}, r) \Vdash \dot{f}(u_{n+1}) = s_{n+1}$.
Let $t_{n+1} \geq r$ such that $\Ht(t_{n+1}) > \Ht(s_{n+1})$.
We leave it to the reader to verify that all of the conditions above hold.

Let $b_0,b_1$ be the downward closure of 
$\{u_n: n \in \omega \}$ and $\{t_{n} : n \in \omega \}$ respectively.
By Lemma \ref{min_extension} there is $q\rq{} < q$ such that
$\alpha_{q'} \geq \delta +1$ and for all branches 
$c \subset T^{q}$, $c$ has an upper bound in $T^{q\rq{}}$ if and only if 
$g_n(b_i)$ is cofinal in $c$ for some  $n \in \omega$ and $ i \in 2$.
Fix such a $q\rq{}$ for the rest of the argument.

We claim that $\Seq{s_n : n \in \omega}$ does not have an 
upper bound in $T^{q'}$.
Suppose for a contradiction that it has an upper bound. 
Then for some $m \in \omega$, either
\begin{enumerate}
\item $\{ g_m (t_n) :n \in \omega \wedge t_n \in \dom(g_m)\}$
is cofinal in the downward closure of $\{s_{n} : n \in \omega \}$ or 
\item $\{ g_m (u_n) : n \in \omega \wedge u_n \in \dom(g_m)\}$ 
is cofinal in the downward closure of $\{s_{n} : n \in \omega \}$.
\end{enumerate}
Due to similarity of the arguments, let\rq{}s assume that the first alternative happens.
Since we enumerated the elements of $\mathcal{F}$ with infinite repetition,
by increasing $m$ if necessary, we can assume that  $t_{m} \in \dom(g_m)$.
But then $g_m(t_{m}) \perp s_{m}$, meaning that
the first alternative cannot happen, which is a contradiction.
Hence $\{s_{n} : n \in \omega \}$ does not have an upper bound in $T^{q\rq{}}$.

Let $t$ be the upper bound of $\Seq{t_{n} : n \in \omega}$ in $T^{q\rq{}}$,
and $u$ be the upper bound  for $\Seq{u_n : n \in \omega}$ which has the lowest height $\delta$.
Let $p$ be a lower bound for $\Seq{p_{n} : n \in \omega}$ which forces that $u \in \dot{S}$.
It is easy to see that
$(q',p,t) \Vdash [\delta \in \dot{C} \wedge u \in \dot{S} \wedge \Ht_{\dot{S}}(u) = \delta].$
Also by  \ref{forcing_relation}, $(q\rq{}, p,t)$ forces 
$\dot{f}(u_n) = s_n$ for all $n \in \omega$.
Hence $(q\rq{},p,t) $ forces that $\dot{f}(u)$ is an upper bound for 
$\Seq{s_{n} : n \in \omega}$ which is a contradiction.
\end{proof}

\def\Dbar{\leavevmode\lower.6ex\hbox to 0pt{\hskip-.23ex \accent"16\hss}D}

\end{document}